\documentclass[12pt,reqno]{amsart}
\usepackage{amsfonts}
\usepackage{amsmath,amsthm,amsfonts,amssymb}
\newtheorem{lemma}{Lemma}[section]
\newtheorem{theorem}{Theorem}[section]
\newtheorem{corollary}{Corollary}[section]

\newtheorem*{thmA}{Theorem A}
\newtheorem*{coroA}{Corollary A}

\def\bl{\begin{lemma}}
\def\bt{\begin{theorem}}
\def\el{\end{lemma}}
\def\et{\end{theorem}}
\def\bp{\begin{proof}}
\def\ep{\end{proof}}
\def\bc{\begin{corollary}}
\def\ec{\end{corollary}}

\def\mb{\mathbb}

\def\O{\Omega}

\def\-{\setminus}

\def\vp{\varphi}

\def\lt{\left}
\def\rt{\right}
\def\+{\bigcup}
\def\.{\bigcap}

\def\ll{\langle}
\def\rl{\rangle}

\title[]
{$L^2$ estimate for polynomials of the Laplace operator with Gaussian measure}

\author []{Shaoyu Dai$^1$, Yang Liu$^2$ and Yifei Pan$^3$}

\address{1 Department of Mathematics, Jinling Institute of Technology, Nanjing, 211169, China.}
\address{\it E-mail address: dymdsy@163.com}

\address{2 Department of Mathematics, Zhejiang Normal University, Jinhua, 321004, China.}
\address{\it E-mail address: liuyang@zjnu.edu.cn}

\address{3 Department of Mathematical Sciences, Purdue University Fort Wayne, Fort Wayne, 46805-1499, USA.}
\address{\it E-mail address: pan@pfw.edu}

\begin{document}
\begin{abstract} Let
$P(\Delta)$
be a polynomial of the Laplace operator $\Delta=\sum_{j=1}^n\frac{\partial^2}{\partial x^2_j}$ on $\mb{R}^n$. We prove
 the existence of weak solutions of the equation $P(\Delta)u=f$ and
the existence of a bounded right inverse of the differential operator $P(\Delta)$ in the weighted Hilbert space with Gaussian measure, i.e., $L^2(\mb{R}^n,e^{-|x|^2})$.
\end{abstract}
\maketitle

\section{Introduction}
In \cite{3}, Rosay gave a simple and elementary proof of the fundamnetal theorem of Malgrange and Ehrenpreis on the existence of fundamental solutions for linear partial differential operators with constant coefficients, and as a by-product, he also proved the following general existence of global solutions for locally square integrable data.
\begin{thmA}
Let $P(D)=\sum a_J\frac{\partial^{|J|}}{\partial x^J}$ be a (nonzero) constant coefficient linear differential operator on $\mb{R}^n$, where $J=(j_1,\cdots,j_n)$, $|J|=j_1+\cdots+j_n$, $\frac{\partial^{|J|}}{\partial x^J}=\frac{\partial^{j_1+\cdots+j_n}}{\partial x_1^{j_1}\cdots\partial x_n^{j_n}}$. Then for every $g\in L^2_{loc}(\mb{R}^n)$ there exists $u\in L^2_{loc}(\mb{R}^n)$ such that $P(D)u=g$.
\end{thmA}

As a corollary of Theorem A, Rosay \cite{3} gives the following result, which means the differential operator $P(D)$ has a bounded right inverse in $L^2(\O)$, where $\O$ is a bounded set in $\mb{R}^n$.
\begin{coroA}
If $\O$ is a bounded set in $\mb{R}^n$, then for every $g\in L^2(\O)$ there exists $u\in L^2(\O)$ such that $P(D)u=g$.
\end{coroA}

However, in this note we will study a bounded right inverse of a special differential operator in an unbounded set.
The special differential operator is
$$P(\Delta)=\Delta^{m}+a_{m-1}\Delta^{m-1}
+\cdots+a_{1}\Delta+a_0,$$
which is a polynomial of the Laplace operator $\Delta=\sum_{j=1}^n\frac{\partial^2}{\partial x^2_j}$ on $\mb{R}^n$, where $m\geq1$ is an integer and $a_{m-1},\cdots,a_{1},a_0$ are complex numbers.
We prove the existence of (global) weak solutions of the equation $P(\Delta)u=f$ in the weighted Hilbert space $L^2(\mb{R}^n,e^{-|x|^2})$ by the following result.

\bt\label{th1.1}
For each $f\in L^2(\mb{R}^n,e^{-|x|^2})$, there exists a weak solution $u\in L^2(\mb{R}^n,e^{-|x|^2})$ solving the equation
$$P(\Delta)u=f$$ in $\mb{R}^n$
with the norm estimate $$\int_{\mb{R}^n} |u|^2e^{-|x|^2}dx\leq
\frac{1}{(8n)^m}\int_{\mb{R}^n} |f|^2e^{-|x|^2}dx.$$
\et

The novelty of Theorem \ref{th1.1} is that the differential operator $P(\Delta)$ has a bounded right inverse
\begin{align*}
Q: L^2(\mb{R}^n,e^{-|x|^2})&\longrightarrow L^2(\mb{R}^n,e^{-|x|^2}), \\
P(\Delta)Q&=I
\end{align*}
with the norm estimate $\|Q\|\leq\frac{1}{(8n)^{\frac{m}{2}}}$.
In particular, the Laplace operator $\Delta$ has a bounded right inverse $Q_0: L^2(\mb{R}^n,e^{-|x|^2})\longrightarrow L^2(\mb{R}^n,e^{-|x|^2})$, which, to the best of our knowledge, appears to be new.

A natural question would be if Theorem \ref{th1.1} would be true for general differential operators as in theorem A. For related results, see \cite{dbar}, \cite{ddd} and \cite{we}.
The method employed in this note was motivated from the H\"{o}rmander $L^2$ method \cite{1} for Cauchy-Riemann
equations from several complex variables.

The organization of this paper is as follows. In Section 2, we will prove several key lemmas based on functional analysis, while the proof of Theorem \ref{th1.1} will be given in Section 3. In Section 4, we will give some further remarks.

\section{Several lemmas}

In this section, we will prove the following key lemma, which is important for the proof of Theorem \ref{th1.1}.
\bl\label{lemma1.1}
Let $\xi$ be a complex number. For each $f\in {L}^2(\mb{R}^n,e^{-|x|^2})$, there exists a weak solution $u\in {L}^2(\mb{R}^n,e^{-|x|^2})$ solving the equation
$$\Delta u+\xi u=f$$ in $\mb{R}^n$
with the norm estimate $$\int_{\mb{R}^n} |u|^2e^{-|x|^2}dx\leq
\frac{1}{8n}\int_{\mb{R}^n} |f|^2e^{-|x|^2}dx.$$
\el

For the proof of Lemma \ref{lemma1.1}, we need some preparation.

First we give some notations.
Here, let ${L}^2_{loc}(\mb{R}^n)$ be the set of all  locally $2$-integrable complex-valued functions on $\mb{R}^n$.
We consider weighted Hilbert space
$${L}^2(\mb{R}^n,e^{-\vp})
=\{f\mid f\in {L}^2_{loc}(\mb{R}^n); \int_{\mb{R}^n}|f|^2 e^{-\vp}dx<+\infty\},$$
where $\vp$ is a nonnegative function on $\mb{R}^n$.
We denote
the weighted inner product for $f,g\in {L}^2(\mb{R}^n,e^{-\vp})$ by
$\langle f,g\rangle_\vp=\int_{\mb{R}^n}\overline{f}g e^{-\vp}dx$
and the weighted norm of $f\in {L}^2(\mb{R}^n,e^{-\vp})$ by
$\|f\|_\vp=\sqrt{\langle f,f\rangle_\vp}.$
Let ${C}_0^\infty(\mb{R}^n)$ denote the set of all smooth complex-valued functions with compact support.
For $u,f\in {L}^2_{loc}(\mb{R}^n)$, we say that $f$ is the Laplace of $u$ in the weak sense, written $\Delta u=f$, provided $\int_{\mb{R}^n}u\Delta\phi dx=\int_{\mb{R}^n}f\phi dx$
for all test functions $\phi\in {C}_0^\infty(\mb{R}^n)$.

In the following, let $\vp$ be a smooth and nonnegative function on $\mb{R}^n$ and $\xi$ be a complex number throughout this section. For $\forall\phi\in {C}_0^\infty(\mb{R}^n)$, we first define the following formal adjoint of $\Delta$ with respect to the weighted inner product in ${L}^2\lt(\mb{R}^n,e^{-\vp}\rt)$. Let $u\in {L}^2_{loc}(\mb{R}^n)$. We calculate as follows.
\begin{align*}
\lt\ll \phi,\Delta u\rt\rl_\vp&=\int_{\mb{R}^n}\overline{\phi} \Delta u e^{-\vp}dx \\
&=\int_{\mb{R}^n}u \Delta\lt(\overline{\phi} e^{-\vp}\rt)dx\\
&=\int_{\mb{R}^n}e^{\vp}u \Delta\lt(\overline{\phi} e^{-\vp}\rt)e^{-\vp}dx\\
&=\lt\ll \overline{e^{\vp}\Delta\lt(\overline{\phi} e^{-\vp}\rt)},u\rt\rl_\vp\\
&=:\lt\ll \Delta_\vp^{*}\phi, u\rt\rl_\vp,
\end{align*}
where $\Delta_\vp^{*}\phi =\overline{e^{\vp}\Delta\lt(\overline{\phi} e^{-\vp}\rt)}$ is the formal adjoint of $\Delta$ with domain in ${C}_0^\infty(\mb{R}^n)$. Let $\lt(\Delta+\xi\rt)_\vp^{*}$ be the formal adjoint of $\Delta+\xi$ with domain in ${C}_0^\infty(\mb{R}^n)$. Note that $I_\vp^{*}=I$, where $I$ is  the identity operator. Then $\lt(\Delta+\xi\rt)_\vp^{*}=\Delta_\vp^{*}+\overline{\xi}$.

Let $\nabla$ be the gradient operator on $\mb{R}^n$.
Now we give several lemmas based on functional analysis.

\bl\label{lemmaifif}
For each $f\in {L}^2(\mb{R}^n,e^{-\vp})$, there exists a global weak solution $u\in {L}^2(\mb{R}^n,e^{-\vp})$ solving the equation
$$\Delta u+\xi u=f$$ in $\mb{R}^n$
with the norm estimate
$$\|u\|^2_\vp\leq c$$
 if and only if
$$|\langle f,\phi\rangle_\vp|^2\leq c\lt\|\lt(\Delta+\xi\rt)^*_\vp\phi\rt\|^2_\vp, \ \ \forall\phi\in {C}_0^\infty(\mb{R}^n),$$
 where $c$ is a constant.
\el

\bp
Let $\Delta+\xi=H$. Then $\lt(\Delta+\xi\rt)^*_\vp=H^*_\vp$.

(Necessity) For $\forall\phi\in {C}_0^\infty(\mb{R}^n)$, from the definition of $H^*_\vp$ and Cauchy-Schwarz inequality, we have
\begin{align*}
|\langle f,\phi\rangle_\vp|^2=|\langle Hu,\phi\rangle_\vp|^2=\lt|\lt\langle u,H^*_\vp\phi\rt\rangle_\vp\rt|^2\leq\|u\|^2_\vp
\lt\|H^*_\vp\phi\rt\|^2_\vp\leq c\lt\|H^*_\vp\phi\rt\|^2_\vp
.
\end{align*}

(Sufficiency) Consider the subspace
$$E=\lt\{H^*_\vp\phi\mid\phi\in {C}_0^\infty(\mb{R}^n)\rt\}\subset {L}^2(\mb{R}^n,e^{-\vp}).$$
Define a linear functional $L_f: E\rightarrow\mb{C}$ by
$$L_f\lt(H^*_\vp\phi\rt)=\langle f,\phi\rangle_\vp=\int_{\mb{R}^n}\overline{f}\phi e^{-\vp}dx.$$
Since
$$\lt|L_f\lt(H^*_\vp\phi\rt)\rt|=\lt|\langle f,\phi\rangle_\vp\rt|
\leq\sqrt{c}\lt\|H^*_\vp\phi\rt\|_\vp,$$
then $L_f$ is a bounded functional on $E$. Let $\overline{E}$ be the closure of $E$ with respect to the norm $\|\cdot\|_\vp$ of ${L}^2(\mb{R}^n,e^{-\vp})$. Note that $\overline{E}$ is a Hilbert subspace of
${L}^2(\mb{R}^n,e^{-\vp})$. So by Hahn-Banach's extension theorem, $L_f$ can be extended to a linear functional $\widehat{L}_f$ on $\overline{E}$
such that
\begin{equation}\label{29}
\lt|\widehat{L}_f(g)\rt|\leq\sqrt{c}\lt\|g\rt\|_\vp, \ \ \forall g\in \overline{E}.
\end{equation}
Using the Riesz representation theorem for $\widehat{L}_f$, there exists a unique $u_0\in \overline{E}$ such that
\begin{equation}\label{30}
\widehat{L}_f(g)=\langle u_0,g\rangle_\vp, \ \ \forall g\in \overline{E}.
\end{equation}

Now we prove $Hu_0=f$. For $\forall\phi\in {C}_0^\infty(\mb{R}^n)$, apply $g=H^*_\vp\phi$ in (\ref{30}). Then
$$\widehat{L}_f\lt(H^*_\vp\phi\rt)=\lt\langle u_0,H^*_\vp\phi\rt\rangle_\vp=\lt\langle Hu_0,\phi\rt\rangle_\vp.$$
Note that
$$\widehat{L}_f\lt(H^*_\vp\phi\rt)=L_f\lt(H^*_\vp\phi\rt)=\langle f,\phi\rangle_\vp.$$
Therefore,
$$\lt\langle Hu_0,\phi\rt\rangle_\vp=\langle f,\phi\rangle_\vp, \ \ \forall\phi\in {C}_0^\infty(\mb{R}^n),$$
i.e.,
$$\int_{\mb{R}^n} \overline{Hu_0}\phi e^{-\vp}dx=\int_{\mb{R}^n} \overline{f}\phi e^{-\vp}dx, \ \ \forall\phi\in {C}_0^\infty(\mb{R}^n).$$
Thus, $Hu_0=f$.

Next we give a bound for the norm of $u_0$. Let $g=u_0$ in (\ref{29}) and (\ref{30}). Then we have
$$\|u_0\|^2_\vp=\lt|\langle u_0,u_0\rangle_\vp\rt|=\lt|\widehat{L}_f(u_0)\rt|
\leq\sqrt{c}\lt\|u_0\rt\|_\vp.$$
Therefore, $\|u_0\|—_\vp^2\leq c$.

Note that $u_0\in\overline{E}$ and $\overline{E}\subset {L}^2(\mb{R}^n,e^{-\vp})$. Then $u_0\in {L}^2(\mb{R}^n,e^{-\vp})$. Let $u=u_0$. So there exists $u\in {L}^2(\mb{R}^n,e^{-\vp})$ such that
$Hu=f$ with $\|u\|^2_\vp\leq c$. The proof is complete.
\ep

\bl\label{lemmaH}
For $\forall\phi\in {C}_0^\infty(\mb{R}^n)$, we have
\begin{align*}
\lt\|\lt(\Delta+\xi\rt)^*_\vp\phi\rt\|^2_\vp
=\lt\|\lt(\Delta+\xi\rt)\phi\rt\|^2_\vp+\lt\langle \phi,\Delta\lt(\Delta_\vp^{*}\phi\rt)-\Delta^{*}_\vp
\lt(\Delta\phi\rt)\rt\rangle_\vp.
\end{align*}
\el

\bp
Let $\Delta+\xi=H$. Then $\lt(\Delta+\xi\rt)^*_\vp=H^*_\vp$.
For $\forall\phi\in {C}_0^\infty(\mb{R}^n)$, we have
\begin{align}
\lt\|H^*_\vp\phi\rt\|^2_\vp&=\lt\langle H^*_\vp\phi,H^*_\vp\phi\rt\rangle_\vp\nonumber\\
&=\lt\langle \phi,HH^*_\vp\phi\rt\rangle_\vp\nonumber\\
&=\lt\langle \phi,H^*_\vp H\phi\rt\rangle_\vp+\lt\langle \phi,HH^*_\vp\phi-H^*_\vp H\phi\rt\rangle_\vp\nonumber\\
&=\lt\langle H\phi,H\phi\rt\rangle_\vp+\lt\langle \phi,HH^*_\vp\phi-H^*_\vp H\phi\rt\rangle_\vp\nonumber\\
&=\lt\|H\phi\rt\|^2_\vp+\lt\langle \phi,HH^*_\vp\phi-H^*_\vp H\phi\rt\rangle_\vp\label{31}
\end{align}
Note that
\begin{align*}
HH^*_\vp\phi&=\lt(\Delta+\xi\rt)\lt(\Delta+\xi\rt)^*_\vp\phi\\
&=\lt(\Delta+\xi\rt)\lt(\Delta^{*}_\vp\phi+\overline{\xi}\phi\rt)\\
&=\Delta\lt(\Delta_\vp^{*}\phi\rt)+\overline{\xi}\Delta\phi+
\xi\Delta_\vp^{*}\phi+|\xi|^2\phi
\end{align*}
and
\begin{align*}
H^*_\vp H\phi&=\lt(\Delta+\xi\rt)^*_\vp\lt(\Delta+\xi\rt)\phi\\
&=\lt(\Delta^{*}_\vp+\overline{\xi}\rt)\lt(\Delta\phi+\xi\phi\rt)\\
&=\Delta^{*}_\vp\lt(\Delta\phi\rt)+\xi\Delta_\vp^{*}\phi+\overline{\xi}\Delta\phi+|\xi|^2\phi.
\end{align*}
Then
\begin{equation}\label{32}
HH^*_\vp\phi-H^*_\vp H\phi=\Delta\lt(\Delta_\vp^{*}\phi\rt)-\Delta^{*}_\vp\lt(\Delta\phi\rt).
\end{equation}
So by (\ref{31}) and (\ref{32}), we have
\begin{align*}
\lt\|H^*_\vp\phi\rt\|^2_\vp=\lt\|H\phi\rt\|^2_\vp+\lt\langle \phi,\Delta\lt(\Delta_\vp^{*}\phi\rt)-\Delta^{*}_\vp\lt(\Delta\phi\rt)\rt\rangle_\vp.
\end{align*}
This lemma is proved.
\ep

\bl\label{lemmak} Let $\vp=|x|^2$. Then for $\forall\phi\in {C}_0^\infty(\mb{R}^n)$, we have
\begin{align*}
\lt\ll\phi, \Delta\lt(\Delta_\vp^{*}\phi\rt)-\Delta^{*}_\vp\lt(\Delta\phi\rt)
\rt\rl_\vp=8n\|\phi\|^2_\vp+8\|\nabla\phi\|^2_\vp.
\end{align*}
\el

\bp
Note that for any smooth functions $\alpha$ and $\beta$ on $\mb{R}^n$, the following formula holds
$$\Delta(\alpha\beta)=\beta\Delta\alpha+\alpha\Delta\beta
+2\nabla\alpha\cdot\nabla\beta.$$
Then for $\forall\phi\in {C}_0^\infty(\mb{R}^n)$, by the definition of $\Delta_\vp^{*}$
we have
\begin{align}\label{1}
\Delta_\vp^{*}\phi=\overline{e^{\vp}\Delta\lt(\overline{\phi} e^{-\vp}\rt)}=\Delta\phi+\phi|\nabla\vp|^2-\phi\Delta\vp-2\nabla\phi\cdot
\nabla\vp.
\end{align}
From (\ref{1})
we have
\begin{align*}
\Delta\lt(\Delta_\vp^{*}\phi\rt)&=\Delta^2\phi+\Delta(\phi|\nabla\vp|^2)
-\Delta(\phi\Delta\vp)-2\Delta(\nabla\phi\cdot
\nabla\vp)\\
&=\Delta^2\phi+\Delta\phi|\nabla\vp|^2+\phi\Delta(|\nabla\vp|^2)+
2\nabla\phi\cdot\nabla(|\nabla\vp|^2)\\
&\ \ \ -\Delta\phi\Delta\vp
-\phi\Delta^2\vp-2\nabla\phi\cdot\nabla(\Delta\vp)
-2\Delta(\nabla\phi\cdot
\nabla\vp)
\end{align*}
and
$$\Delta^*_\vp(\Delta\phi)=\Delta^2\phi+\Delta\phi|\nabla\vp|^2
-\Delta\phi\Delta\vp-2\nabla(\Delta\phi)\cdot\nabla\vp.$$
Then
\begin{align}
\Delta\lt(\Delta_\vp^{*}\phi\rt)-\Delta^{*}_\vp\lt(\Delta\phi\rt)&=
\phi\Delta(|\nabla\vp|^2)+
2\nabla\phi\cdot\nabla(|\nabla\vp|^2)
-\phi\Delta^2\vp\nonumber\\&\ \ \ -2\nabla\phi\cdot\nabla(\Delta\vp)
-2\Delta(\nabla\phi\cdot
\nabla\vp)+2\nabla(\Delta\phi)\cdot\nabla\vp\label{41}.
\end{align}
Let $\vp=|x|^2$. We have $\nabla\vp=2x$, $\Delta\vp=2n$,
$|\nabla\vp|^2=4|x|^2$, $\nabla(|\nabla\vp|^2)=8x$, $\Delta(|\nabla\vp|^2)=8n$.
Then by (\ref{41}) and the following formula
$$\Delta(\nabla\phi\cdot x)=\nabla(\Delta\phi)\cdot x+2\Delta\phi,\ \ \forall\phi\in {C}_0^\infty(\mb{R}^n),$$ we get
\begin{align*}
\Delta\lt(\Delta_\vp^{*}\phi\rt)-\Delta^{*}_\vp\lt(\Delta\phi\rt)
=8n\phi+16(\nabla\phi\cdot x)-8\Delta\phi.
\end{align*}
Consequently,
\begin{align*}
\lt\ll\phi, \Delta\lt(\Delta_\vp^{*}\phi\rt)-\Delta^{*}_\vp\lt(\Delta\phi\rt)
\rt\rl_\vp&=\lt\ll\phi, 8n\phi+16(\nabla\phi\cdot x)-8\Delta\phi
\rt\rl_\vp\\
&=8n\|\phi\|^2_\vp+8\lt\ll\phi, 2(\nabla\phi\cdot x)-\Delta\phi
\rt\rl_\vp.
\end{align*}
Note, as the key step of the proof, that
\begin{align*}
\lt\ll\phi, 2(\nabla\phi\cdot x)-\Delta\phi
\rt\rl_\vp&=\int_{\mb{R}^n}\overline{\phi}(2(\nabla\phi\cdot x)-\Delta\phi)e^{-\vp}dx\\
&=\int_{\mb{R}^n}\overline{\phi}\sum^n_{j=1}
\lt(2x_j\frac{\partial\phi}{\partial x_j}-\frac{\partial^2\phi}{\partial x_j\partial x_j}\rt)e^{-|x|^2}dx\\
&=-\int_{\mb{R}^n}\overline{\phi}\sum^n_{j=1}\frac{\partial}{\partial x_j}\lt(\frac{\partial\phi}{\partial x_j}e^{-|x|^2}\rt)dx\\
&=-\sum^n_{j=1}\int_{\mb{R}^n}\overline{\phi}\frac{\partial}{\partial x_j}\lt(\frac{\partial\phi}{\partial x_j}e^{-|x|^2}\rt)dx\\
&=\sum^n_{j=1}\int_{\mb{R}^n}\frac{\partial\overline{\phi}}{\partial x_j}\lt(\frac{\partial\phi}{\partial x_j}e^{-|x|^2}\rt)dx\\
&=\sum^n_{j=1}\int_{\mb{R}^n}\lt|\frac{\partial\phi}{\partial x_j}\rt|^2e^{-|x|^2}dx\\
&=\int_{\mb{R}^n}|\nabla\phi|^2e^{-|x|^2}dx\\
&=\|\nabla\phi\|^2_\vp.
\end{align*}
Then
\begin{align*}
\lt\ll\phi, \Delta\lt(\Delta_\vp^{*}\phi\rt)-\Delta^{*}_\vp\lt(\Delta\phi\rt)
\rt\rl_\vp=8n\|\phi\|^2_\vp+8\|\nabla\phi\|^2_\vp.
\end{align*}
The lemma is proved.
\ep

Now we give the proof of Lemma \ref{lemma1.1}.

\bp
Let $\vp=|x|^2$. By Lemma \ref{lemmaH} and Lemma \ref{lemmak}, we have for $\forall\phi\in {C}_0^\infty(\mb{R}^n)$,
\begin{align}
\lt\|\lt(\Delta+\xi\rt)^*_\vp\phi\rt\|^2_\vp
&=\lt\|\lt(\Delta+\xi\rt)\phi\rt\|^2_\vp+\lt\langle \phi,\Delta\lt(\Delta_\vp^{*}\phi\rt)-\Delta^{*}_\vp
\lt(\Delta\phi\rt)\rt\rangle_\vp\nonumber\\
&\geq\lt\langle \phi,\Delta\lt(\Delta_\vp^{*}\phi\rt)-\Delta^{*}_\vp
\lt(\Delta\phi\rt)\rt\rangle_\vp\nonumber\\
&\geq8n\|\phi\|^2_\vp.\label{34}
\end{align}
By Cauchy-Schwarz inequality and (\ref{34}), we have for $\forall\phi\in {C}_0^\infty(\mb{R}^n)$,
\begin{align*}
|\langle f,\phi\rangle_\vp|^2
&\leq\lt\|f\rt\|^2_\vp
\lt\|\phi\rt\|^2_\vp\\
&=\lt(\frac{1}{8n}\lt\|f\rt\|^2_\vp\rt)
\lt(8n\lt\|\phi\rt\|^2_\vp\rt)\\
&\leq\lt(\frac{1}{8n}\lt\|f\rt\|^2_\vp\rt)
\lt\|\lt(\Delta+\xi\rt)^*_\vp\phi\rt\|^2_\vp.
\end{align*}
Let $c=\frac{1}{8n}\lt\|f\rt\|^2_\vp$. Then
$$|\langle f,\phi\rangle_\vp|^2\leq c\lt\|\lt(\Delta+\xi\rt)^*_\vp\phi\rt\|^2_\vp,\ \ \forall\phi\in {C}_0^\infty(\mb{R}^n).$$
By Lemma \ref{lemmaifif}, there exists a global weak solution $u\in {L}^2(\mb{R}^n,e^{-\vp})$ solving the equation
$$\Delta u+\xi u=f$$ in $\mb{R}^n$
with the norm estimate
$$\|u\|^2_\vp\leq c,$$
i.e.,
$$\Delta u+\xi u=f \ \ {\it with} \ \ \int_{\mb{R}^n} |u|^2e^{-|x|^2}dx\leq
\frac{1}{8n}\int_{\mb{R}^n} |f|^2e^{-|x|^2}dx.$$
The proof is complete.
\ep

\section{Proof of Theorem \ref{th1.1}}

Now we give the proof of Theorem \ref{th1.1}.
\bp
Let $\vp=|x|^2$. By the fundamental theorem of algebra, the polynomial $P(\Delta)$ can be rewritten as
$$P(\Delta)=\lt(\Delta+\xi_1\rt)\cdots\lt(\Delta+\xi_m\rt),$$
where $\xi_j$ is a complex number for $j=1,\cdots,m$.

Let $H_j=\Delta+\xi_j$ for $j=1,\cdots,m$. If $m=1$, then the theorem is proved by Lemma \ref{lemma1.1}. Now assume that $m\geq2$. For $f$ and $H_1$, by Lemma \ref{lemma1.1}, there exists $u_1\in {L}^2(\mb{R}^n,e^{-|x|^2})$ such that
\begin{align*}
H_1u_1=f \ \ {\it with} \ \ \int_{\mb{R}^n} |u_1|^2e^{-|x|^2}dx\leq
\frac{1}{8n}\int_{\mb{R}^n} |f|^2e^{-|x|^2}dx.
\end{align*}
For $u_1$ and $H_2$, by Lemma \ref{lemma1.1}, there exists $u_2\in {L}^2(\mb{R}^n,e^{-|x|^2})$ such that
\begin{align*}
H_2u_2=u_1 \ \ {\it with} \ \ \int_{\mb{R}^n} |u_2|^2e^{-|x|^2}dx\leq
\frac{1}{8n}\int_{\mb{R}^n} |u_1|^2e^{-|x|^2}dx.
\end{align*}
So by the same method, we have for $1\leq j\leq m-1$, $u_j$ and $H_{j+1}$,
there exists $u_{j+1}\in {L}^2(\mb{R}^n,e^{-|x|^2})$ such that
\begin{align*}
H_{j+1}u_{j+1}=u_j \ \ {\it with} \ \ \int_{\mb{R}^n} |u_{j+1}|^2e^{-|x|^2}dx\leq
\frac{1}{8n}\int_{\mb{R}^n} |u_j|^2e^{-|x|^2}dx.
\end{align*}
Thus, there exists $u_m\in {L}^2(\mb{R}^n,e^{-|x|^2})$ such that
\begin{align*}
H_1\cdots H_mu_m=f \ \ {\it with} \ \ \int_{\mb{R}^n} |u_m|^2e^{-|x|^2}dx\leq
\frac{1}{(8n)^m}\int_{\mb{R}^n} |f|^2e^{-|x|^2}dx,
\end{align*}
i.e.,
\begin{align*}
P(\Delta)u_m=f \ \ {\it with} \ \ \int_{\mb{R}^n} |u_m|^2e^{-|x|^2}dx\leq
\frac{1}{(8n)^m}\int_{\mb{R}^n} |f|^2e^{-|x|^2}dx,
\end{align*}

Let $u_m=u$. Then by the above formula, the theorem is proved.
\ep

\section{Further remarks}

\noindent\textbf{Remark 1.}
Given $\lambda>0$ and $x_0\in\mb{R}^n$, for the weight $\vp=\lambda|x-x_0|^2$, we obtain the following corollary from Theorem \ref{th1.1}.

\bc\label{thzy}
For each $f\in L^2(\mb{R}^n,e^{-\lambda|x-x_0|^2})$, there exists a weak solution
$u\in L^2(\mb{R}^n,e^{-\lambda|x-x_0|^2})$ solving the equation
$$P(\Delta)u=f$$
with the norm estimate $$\int_{\mb{R}^n} |u|^2e^{-\lambda|x-x_0|^2}dx\leq
\frac{1}{\lambda^{2m}(8n)^m}\int_{\mb{R}^n} |f|^2e^{-\lambda|x-x_0|^2}dx.$$
\ec

\bp
From $f\in L^2(\mb{R}^n,e^{-\lambda|x-x_0|^2})$, we have
\begin{equation}\label{101}
\int_{\mb{R}^n} |f|^2(x)e^{-\lambda|x-x_0|^2}dx<+\infty.
\end{equation}
Let $x=\frac{y}{\sqrt{\lambda}}+x_0$ and $g(y)=f(x)= f\lt(\frac{y}{\sqrt{\lambda}}+x_0\rt)$. Then by (\ref{101}), we have
\begin{equation*}
\frac{1}{\lt(\sqrt{\lambda}\rt)^n}\int_{\mb{R}^n} |g(y)|^2 e^{-|y|^2}dy<+\infty,
\end{equation*}
which implies that $g\in L^2(\mb{R}^n,e^{-|y|^2})$. For $g$, applying Theorem \ref{th1.1} with $P(\Delta)$ replaced by
$$\widetilde{P}(\Delta)=\Delta^{m}+\frac{a_{m-1}}{\lambda}\Delta^{m-1}+\frac{a_{m-2}}{\lambda^2}\Delta^{m-2}
+\cdots+\frac{a_{1}}{\lambda^{m-1}}\Delta+\frac{a_0}{\lambda^m}$$
there exists a weak solution $v\in L^2(\mb{R}^n,e^{-|y|^2})$ solving the equation
\begin{equation}\label{103}
\widetilde{P}(\Delta)v(y)=g(y)
\end{equation}
in $\mb{R}^n$ with the norm estimate
\begin{equation}\label{104}
\int_{\mb{R}^n} |v(y)|^2e^{-|y|^2}dy\leq
\frac{1}{(8n)^m}\int_{\mb{R}^n} |g(y)|^2e^{-|y|^2}dy.
\end{equation}
Note that $y=\sqrt{\lambda}(x-x_0)$ and $g(y)=f(x)$. Let $u(x)=\frac{1}{\lambda^m}v(y)=\frac{1}{\lambda^m}v\lt(\sqrt{\lambda}(x-x_0)\rt)$. Then (\ref{103}) and (\ref{104}) can be rewritten by
\begin{equation}\label{109}
P(\Delta)u(x)
=f(x)
\end{equation}
\begin{equation}\label{110}
\int_{\mb{R}^n} |u(x)|^2 e^{-\lambda|x-x_0|^2}dx\leq
\frac{1}{\lambda^{2m}(8n)^m}\int_{\mb{R}^n} |f(x)|^2e^{-\lambda|x-x_0|^2}dx.
\end{equation}
(\ref{110}) implies that $u\in L^2(\mb{R}^n,e^{-\lambda|x-x_0|^2})$. Then by (\ref{109}) and (\ref{110}), the proof is complete.
\ep

\noindent\textbf{Remark 2.} When $f\in L^2(\mb{R}^n)$, the solutins of $P(\Delta)u=f$ are not necessary in $L^2(\mb{R}^n)$. For example: $n=1, P(\Delta)=\Delta$,
\begin{align*}
f(x)=\left\{
\begin{array}{ccc}
\frac{1}{x},       &      & {x\geq1,}\\
x,     &      & {0<x<1,}\\
0,     &      & {x\leq0,}
\end{array} \right.
\end{align*}
$$u(x)=\int^x_0(x-t)f(t)dt+c_1x+c_2
=-\frac{x}{2}+xlnx+\frac{2}{3}+c_1x+c_2,\ \ x\geq1,$$
where $c_1$ and $c_2$ are arbitrary real constants. It is easy to see $u\not\in L^2(\mb{R})$.

\end{document}